\documentstyle[11pt,amstex,amssymb]{amsart}\textheight=19.8truecm\textwidth=13.7truecm
\begin{document} 

\vskip 30pt
\centerline{\Large\bf Hermitian operators and convex functions}
\vskip 20pt
\centerline{Jean-Christophe Bourin}
\vskip 10pt
 \centerline{E-mail: bourinjc@@club-internet.fr}
 \vskip 5pt
 \centerline{Universit\'e de Cergy-Pontoise, d\'ept.\ de Math\'ematiques}
 \vskip 5pt
 \centerline{2 rue Adolphe Chauvin, 95302 Pontoise, France}

\vskip 15pt
\noindent
{\small {\bf Abstract.} 
\vskip 5pt
We establish several convexity results for Hermitian matrices.  For instance: Let  $A$, $B$ be Hermitian  and let $f$ be a convex function. If $X$ and $Y$ stands for $f(\{A+B\}/2)$ and $\{f(A)+f(B)\}/2$ respectively, then there exist unitaries $U$, $V$ such that
$$
X \le \frac{UYU^* + VYV^*}{2}.
$$
This is nothing but the matrix version of the  scalar convexity inequality 
$$
f\left(\frac{a+b}{2}\right) \le \frac{f(a)+f(b)}2.
$$

As a consequence we get,  \, $\lambda_{2j-1}(X) \le \lambda_j(Y)$,\, where $\lambda_j(\cdot)$ are the eigenvalues arranged in decreasing order. 


Mathematical subjects classification:   47A30 47A63}

\vskip 15pt
{\large\bf Introduction}
\vskip 10pt
The main aim of this paper is to give a matrix version of the scalar inequality
\begin{equation}
f\left(\frac{a+b}{2}\right) \le \frac{f(a)+f(b)}2
\end{equation}
for convex functions $f$ on the real line.

Capital letters $A$, $B\dots Z$ mean $n$-by-$n$ complex matrices, or operators on a finite dimensional Hilbert space ${\cal H}$; $I$ stands for the identity. When $A$ is positive semidefinite, resp.\ positive definite, we write $A\ge 0$, resp.\ $A>0$. 

A classical matrix version of (1) is von Neuman's Trace Inequality: For Hermitians $A$, $B$,
\begin{equation}
{\rm Tr}\,f\left(\frac{A+B}{2}\right) \le {\rm Tr}\,\frac{f(A)+f(B)}{2}
\end{equation}
When $f$ is convex and monotone, we showed [2] that  (2) can be extended to an operator inequality: There exists a unitary $U$ such that
\begin{equation}
f\left(\frac{A+B}{2}\right) \le U\cdot\frac{f(A)+f(B)}{2}\cdot U^*
\end{equation}
We also established  similar inequalities involving more general convex combinations. These inequalities are  equivalent to an inequality for compressions. Recall that given an operator $Z$ and a subspace ${\cal E}$ with corresponding orthoprojection $E$, the compression of $Z$ onto ${\cal E}$, denoted by $Z_{\cal E}$, is the restriction of $EZ$ to ${\cal E}$. Inequality (3) can be derived from: For every 
Hermitian $A$, subspace ${\cal E}$ and monotone convex function $f$, there exists a unitary operator $U$ on ${\cal E}$ such that
\begin{equation}
f(A_{\cal E}) \le Uf(A)_{\cal E} U^*.
\end{equation}
Inequalities (3) and (4) are equivalent to inequalities for eigenvalues. For instance (4) can be rephrased as
\begin{equation*}
\lambda_j(f(A_{\cal E})) \le \lambda_j(f(A)_{\cal E}), \quad \ j=1,\, 2,\dots
\end{equation*}
where $\lambda_j(\cdot)$, $j=1,\,2,\dots$ are the eigenvalues arranged in decreasing order and counted with their multiplicities.
Having proved an inequality such as (3) for monotone convex functions, it remains to search counterparts for general convex functions. We derived from (3) the following result for even convex functions $f$\,: Given Hermitians $A$,\,$B$, there exist unitaries $U$, $V$ such that
\begin{equation}
f\left(\frac{A+B}{2}\right) \le \,\frac{Uf(A)U^*+Vf(B)V^*}{2}.
\end{equation}
This generalizes a wellknown inequality for the absolute value,
$$
|A+B|\le U|A|U^*+V|B|V^*
$$
We do not know whether (5) is valid for all convex functions. 

In Section 1 we present a counterpart of (4) for all convex functions. This will enable us to give, in Section 2, a quite natural counterpart of (3) for all convex functions. Though (3) can be proven independently of (4) -and the same for the counterparts-, we have the feeling that in the case of general convex functions, the approach via compressions is more illuminating.

\vskip 20 pt\noindent
{\large\bf 1. Compressions}
\vskip 10pt

Our substitute to (4) for general convex functions (on the real line) is:

\vskip 10pt \noindent
{\bf Theorem 1.1.} {\it Let  $A$ be  Hermitian, let   
 ${\cal E}$ be a subspace and let $f$ be  a convex function. Then, there exist unitaries  $U$, $V$ on ${\cal E}$ such that
\begin{equation*}
f(A_{\cal E}) \le \frac{Uf(A)_{\cal E}U^* + Vf(A)_{\cal E}V^*}{2}.
\end{equation*}
Consequently, for $j=1,\,2,\dots$,
$$
\lambda_{2j-1}(f(A_{\cal E})) \le \lambda_j(f(A)_{\cal E}).
$$
}

\newpage

\vskip 10pt\noindent
{\bf Proof.} We may find spectral subspaces ${\cal E}'$ and ${\cal E}''$ for $A_{\cal E}$ and a real $r$ such that

 \ \ (i) ${\cal E}={\cal E}'\oplus{\cal E}''$,
 
\ (ii) the spectrum of $A_{\cal E'}$ lies on $(-\infty,r]$ and the spectrum of  $A_{\cal E''}$ lies on $[r,\infty)$,
 
 (iii) $f$ is monotone both on $(-\infty,r]$ and  $[r,\infty)$.

Let $k$ be an integer, $1\le k\le \dim{\cal E}'$.  There exists a spectral subspace ${\cal F}\subset{\cal E'}$ for $A_{\cal E'}$ (hence for $f(A_{\cal E'})$), $\dim {\cal F}=k$, such that
\begin{align*} \lambda_k[f(A_{\cal E'})] &=\min_{h\in{\cal F};\ \Vert h\Vert=1} \langle h,f(A_{\cal F})h \rangle  \\
&= \min\{f(\lambda_1(A_{\cal F}))\,;\,f(\lambda_k(A_{\cal F}))\} \\
&= \min_{h\in{\cal F};\ \Vert h\Vert=1} f(\langle h,A_{\cal F}h \rangle)  \\
&= \min_{h\in{\cal F};\ \Vert h\Vert=1} f(\langle h,Ah \rangle) 
\end{align*}
where at the second and third steps we use the monotony of $f$ on $(-\infty,r]$ and the fact that $A_{\cal F}$'s spectrum lies on $(-\infty,r]$. The convexity of $f$ implies 
$$
f(\langle h,Ah \rangle) \le \langle h,f(A)h \rangle
$$
for all normalized vectors $h$. Therefore, by the minmax principle,
\begin{align*}
\lambda_k[f(A_{\cal E'})] &\le \min_{h\in{\cal F};\ \Vert h\Vert=1} \langle h,f(A)h \rangle \\
&\le \lambda_k[f(A)_{\cal E'}].
\end{align*} 
This  statement is equivalent (by unitay congruence to diagonal matrices) to the existence of a unitary operator $U_0$ on ${\cal E}'$ such that 
$$
f(A_{\cal E'})\le U_0 f(A)_{\cal E'}U_0^*.
$$
Similarly we get a unitary $V_0$ on ${\cal E''}$ such that
$$
f(A_{\cal E''})\le V_0 f(A)_{\cal E''}V_0^*.
$$
Thus we have
$$
f(A_{\cal E})\le
\begin{pmatrix} 
U_0 &0 \\ 0&V_0 
\end{pmatrix}  
 \begin{pmatrix} 
f(A)_{\cal E'} &0 \\ 0& f(A)_{\cal E''}
\end{pmatrix}  
 \begin{pmatrix} 
U_0^* &0 \\ 0&V_0^* 
\end{pmatrix}. 
$$
Besides we note that, still in respect with the decomposition ${\cal E}={\cal E}'\oplus{\cal E}''$,
$$
\begin{pmatrix} 
f(A)_{\cal E'} &0 \\ 0& f(A)_{\cal E''} 
\end{pmatrix}
= \frac{1}{2}\left\{  
\begin{pmatrix} 
I &0 \\ 0& I 
\end{pmatrix}
f(A)_{\cal E}
\begin{pmatrix} 
I &0 \\ 0& I 
\end{pmatrix}
+
\begin{pmatrix} 
I &0 \\ 0& -I 
\end{pmatrix}
f(A)_{\cal E}
\begin{pmatrix} 
I &0 \\ 0& -I 
\end{pmatrix}
\right\}.
$$ 
So, letting
$$
U=\begin{pmatrix} 
U_0 &0 \\ 0& V_0 
\end{pmatrix}
\quad{\rm and}\quad
V=
\begin{pmatrix} 
U_0 &0 \\ 0& -V_0 
\end{pmatrix}
$$
we get
\begin{equation}
f(A_{\cal E}) \le \frac{Uf(A)_{\cal E}U^* + Vf(A)_{\cal E}V^*}{2}.
\end{equation}

It remains to check that (6) entails
$$
\lambda_{2j-1}(f(A_{\cal E})) \le \lambda_j(f(A)_{\cal E}).
$$
This follows from the forthcoming elementary observation. \qquad $\Box$

\vskip 10pt \noindent
{\bf Proposition 1.2.} {\it Let  $X$, $Y$ be  Hermitians such that  
\begin{equation}
X \le \frac{UYU^* + VYV^*}{2}
\end{equation}
for some unitaries $U$, $V$. Then, for $j=1,\,2,\dots$,
$$
\lambda_{2j-1}(X) \le \lambda_j(Y).
$$
}

\vskip 10pt \noindent
{\bf Proof.} By adding a $rI$ term, for a suitable scalar $r$, both to $X$ and $Y$, it suffices to show that
\begin{equation}
\lambda_{2j-1}(X)>0 \quad \Longrightarrow \quad \lambda_j(Y)>0.
\end{equation}
We need the following obvious fact: Given Hermitians $A$, $B$,
$$
{\rm rank} (A+B)_+ \le {\rm rank} A_+ +{\rm rank} B_+ 
$$
where the subscript + stands for positive parts. Applying this to $A=UYU^*$ and $B=VYV^*$ we infer that the negation of (8), that is $\lambda_{2j-1}(A+B)>0$ and $\lambda_j(A)\,(=\lambda_j(B))\le 0$, can not hold. Indeed, the relation
$$
{\rm rank} (A+B)_+ \ge 2j-1 > (j-1)+(j-1) \ge {\rm rank} A_+ +{\rm rank} B_+ 
$$
would contradict the previous rank inequality. \qquad $\Box$

\vskip 10pt \noindent
{\bf Remark 1.3.} From inequality (7) one also derives, as a straightforward consequence of Fan's Maximum Principle [1, Chapter 4], 
$$
\sum_{j=1}^k\lambda_j(X) \le \sum_{j=1}^k\lambda_j(Y)
$$
for $k=1,\,2,\dots$.

 Inequality (7) also implies
$$
\lambda_{i+j+1}(X) \le \frac{1}{2}\{\lambda_{i+1}(Y) +\lambda_{j+1}(Y)\}
$$
for $i,\,j=0,\,1,\,\dots$. It is a special case of  Weyl's inequalities [1, Chapter 3].
\vskip 10pt \noindent
{\bf Remark 1.4.} For operators acting on an infinite dimensional (separable) space, the main inequality of Theorem 1.1 is still valid at the cost of an additional $rI$ term in the RHS, with $r>0$ arbitrarily small. See [3, Chapter 1] for the analogous result for (4).

\vskip 10pt
Obviously, for a concave function $f$, the main inequality of Theorem 1.1 is reversed. But the following is open:

\vskip 10pt \noindent
{\bf Question 1.5.} Let $g$ be a concave function, let $A$ be Hermitian and let ${\cal E}$ be a subspace.
Can we find unitaries $U$, $V$ on ${\cal E}$ such that
\begin{equation*}
g(A)_{\cal E} \le \frac{Ug(A_{\cal E})U^* + Vg(A_{\cal E})V^*}{2} \ \ ?
\end{equation*}

\vskip 20 pt\noindent
{\large\bf 2. Convex combinations}
\vskip 10pt
The next two theorems can be regarded as matrix versions of Jensen's inequality. The first one is also a matrix version of the elementary scalar inequality
$$
f(za)\le zf(a)
$$
for convex functions $f$ with $f(0)\le 0$ and scalars $a$ and  $z$ with $0<z<1$.

\vskip 10pt \noindent
{\bf Theorem 2.1.} {\it Let $f$ be a convex function, let $A$ be  Hermitian, let $Z$ be a contraction and set $X=f(Z^*AZ)$ and $Y= Z^*f(A)Z$. Then, there exist  unitaries  $U$, $V$ such that
\begin{equation*}
X \le \frac{UYU^* + VYV^*}{2}.
\end{equation*}
}

\vskip 10pt
A family $\{Z_i\}_{i=1}^m$ is an isometric column if $\sum_{i=1}^m Z_i^*Z_i=I$

\vskip 10pt \noindent
{\bf Theorem 2.2.} {\it Let $f$ be a convex function, let $\{A_i\}_{i=1}^m$ be  Hermitians, let $\{Z_i\}_{i=1}^m$ be an isometric column and set $X=f(\sum Z_i^*A_iZ_i)$ and $Y=\sum Z_i^*f(A_i)Z_i$. Then, there exist  unitaries  $U$, $V$ such that
\begin{equation*}
X \le \frac{UYU^* + VYV^*}{2}.
\end{equation*}
}

\vskip 10pt \noindent
{\bf Corollary 2.3.} {\it Let $f$ be a convex function, let $A$, $B$ be  Hermitians and  set $X=f(\{A+B\}/2)$ and $Y=\{f(A)+f(B)\}/2$. Then, there exist  unitaries  $U$, $V$ such that
\begin{equation*}
X \le \frac{UYU^* + VYV^*}{2}.
\end{equation*}
}

\vskip 10pt\noindent
Recall that the above inequality entails that for $j=1,\,2,\dots$,
 $$\lambda_{2j-1}(X) \le \lambda_j(Y).$$

\vskip 10pt
 We turn to the proof of Theorems 2.1 and 2.2.
\vskip 10pt \noindent
{\bf Proof.} Theorem 2.1 and  Theorem 1.1 are equivalent. Indeed, to prove Theorem 1.1,  we may assume  that  $f(0)=0$. Then, Theorem 1.1 follows from Theorem 2.2 by taking $Z$ as the projection onto ${\cal E}$.

 Theorem 1.1 entails Theorem 2.2: to see that, we introduce the partial isometry $J$ and the operator $\tilde{A}$ on ${\cal H}\oplus{\cal H}$ defined by
$$ J=
\begin{pmatrix}
Z&0 \\ (I-|Z|^2)^{1/2}&0
\end{pmatrix},\quad
\tilde{A}=
\begin{pmatrix}
A&0 \\ 0&0
\end{pmatrix}.
$$
Denoting by ${\cal H}$ the first summand of the direct sum ${\cal H}\oplus{\cal H}$, we observe that
$$
f(Z^*AZ)= f(J^*\tilde{A}J)\!:\!{\cal H} = J^*f(\tilde{A}_{J({\cal H})})J\!:\!{\cal H}.
$$
where $X\!:\!{\cal H}$ means the restriction of an operator $X$ to the first summand of ${\cal H}\oplus{\cal H}$.
Applying  Theorem 1.1 with ${\cal E}=J({\cal H})$, we get  unitaries  $U_0$, $V_0$ on  $J({\cal H})$ such that
$$
f(Z^*AZ) \le J^*\frac{ U_0f(\tilde{A})_{J({\cal H})}U_0^* +
V_0f(\tilde{A})_{J({\cal H})}V_0^*}{2}J \!:\!{\cal H}. 
$$
Equivalently, there exist  unitaries  $U$, $V$ on ${\cal H}$ such that
\begin{align*}
f(Z^*AZ) \le& \frac{UJ^*f(\tilde{A})_{J({\cal H})}(J\!:\!{\cal H})U^*
+VJ^*f(\tilde{A})_{J({\cal H})}(J\!:\!{\cal H})V^*}{2}
 \\
=&\frac{1}{2}\left\{UJ^* \begin{pmatrix} f(A)&0 \\0&f(0) \end{pmatrix}(J\!:\!{\cal H})U^*+
VJ^* \begin{pmatrix} f(A)&0 \\0&f(0) \end{pmatrix}(J\!:\!{\cal H})V^*\right\} \\
=& \frac{1}{2}U\{Z^*f(A)Z + (I-|Z|^2)^{1/2}f(0)(I-|Z|^2)^{1/2}\}U^* \\
&\qquad\qquad\qquad + \frac{1}{2}V\{Z^*f(A)Z + (I-|Z|^2)^{1/2}f(0)(I-|Z|^2)^{1/2}\}V^*.
\end{align*}
Using $f(0)\le 0$ we obtain the first claim of Theorem 2.2.

 Similarly, Theorem 1.1 implies Theorem 2.3 (we may assume $f(0)=0$) by considering the partial isometry and the operator on $\oplus^m{\cal H}$,
$$
\begin{pmatrix}
Z_1&0&\cdots&0 \\ \vdots &\vdots &\ &\vdots \\
Z_m &0 &\cdots &0
\end{pmatrix},\qquad
\begin{pmatrix}
A_1&\ &\  \\ 
\ &\ddots &\  \\
\ &\  &A_m
\end{pmatrix} .
$$
\qquad $\Box$

\vskip 10pt
We note that our theorems contain two wellknown trace inequalities [4], [5]:

\vskip 10pt\noindent
{\bf 2.4. Brown-Kosaki:} Let $f$ be convex with  $f(0)\le 0$ and let $A$ be Hermitian. Then, for  all contractions $Z$, 
\begin{equation*}
{\rm Tr}\,f(Z^*AZ) \le {\rm Tr}\,Z^*f(A)Z. 
\end{equation*}

\vskip 10pt\noindent
 {\bf 2.5. Hansen-Pedersen:}   Let $f$ be convex  and let $\{A_i\}_{i=1}^m$ be Hermitians.  Then, for  all isometric column $\{Z_i\}_{i=1}^m$,
\begin{equation*}
{\rm Tr}\,f(\sum_iZ_i^*A_iZ_i) \le {\rm Tr}\,\sum_iZ_i^*f(A_i)Z_i. 
\end{equation*}

\vskip 30pt
{\bf References}

\noindent
{\small 

\noindent
[1] R.\ Bhatia, Matrix Analysis, Springer, Germany, 1996.

\noindent
[2] J.-C.\ Bourin, {\it Convexity or concavity inequalities for Hermitian operators},   Math.\ Ineq.\ Appl. 7 ${\rm n}^04$  (2004) 607-620.

\noindent
[3] J.-C.\ Bourin, Compressions, Dilations and Matrix Inequalities, RGMIA monograph, Victoria university, 
Melbourne 2004 (http://rgmia.vu.edu.au/monograph)

\noindent
[4] L.\ G.\ Brown and H.\ Kosaki, {\it Jensen's inequality in semi-finite von Neuman algebras}, J. Operator Theory 23 (1990) 3-19.

\noindent
[5] F.\ Hansen G.\ K.\ Pedersen, {\it Jensen's operator inequality}, Bull. London Math.\ Soc.\
35 (2003) 553-564.

\end{document}